\documentstyle[12pt]{article}
\title{ The  Radicals of Crossed Products
 \thanks {This work was supported by the  National Natural Science Foundation  }}
\author{
Shouchuan Zhang  \\ Department  of Mathematics, Hunan University\\
Changsha  410082, \
 P.R.China. \ \
E-mail:z9491@yahoo.com.cn\\
}
\date{}
\begin{document}
\newtheorem{Theorem}{\quad Theorem}[section]
\newtheorem{Proposition}[Theorem]{\quad Proposition}
\newtheorem{Definition}[Theorem]{\quad Definition}
\newtheorem{Corollary}[Theorem]{\quad Corollary}
\newtheorem{Lemma}[Theorem]{\quad Lemma}
\newtheorem{Example}[Theorem]{\quad Example}
\maketitle \addtocounter{section}{-1}

 \begin {abstract} The relations between the radical of crossed
product $R \# _\sigma H$ and algebra $R$ are obtained. Using this
theory, the author shows that if $H$ is a finite-dimensional
semisimple, cosemisimle, and either commutative
   or cocommutative Hopf algebra, then $R$ is $H$-semiprime iff
      $R$ is semiprime iff
$R\#_{\sigma}H$ is semiprime.

 \end {abstract}

\section {Introduction and Preliminaries}

       J.R. Fisher \cite{Fi75} built up the general theory of $H$-radicals for
$H$-module algebras. He studied $H$-Jacobson radical and  obtained
\begin {eqnarray}
  r_{j}(R\#H) \cap R  = r_{Hj}(R)
 \label {e (1)}
 \end {eqnarray}
\noindent for any irreducible Hopf algebra $H$(\cite [Theorem
4]{Fi75}).
 J.R. Fisher \cite{Fi75} asked when is
 \begin {eqnarray}
 r_{j}(R\#H) =r_{Hj}(R) \#H
\label {e (2)}
\end {eqnarray}
\noindent and asked if
\begin {eqnarray}
 r_{j}(R\#H) \subseteq (r_{j}(R):H) \#H
\label {e (3)}
\end {eqnarray}
R.J. Blattner, M. Cohen and S. Montgomery in \cite {BCM86} asked
whether  $R\#_{\sigma}H$ is semiprime with a finite-dimensional
 semisimple Hopf algebra $H$ when $R$ is semiprime,
 which is called the semiprime problem.

   If $H$ is a finite-dimensional semisimple Hopf algebra and $R$ is semiprime,
    then $R\#_{\sigma}H$ is semiprime
in  the following five cases:

(i)  $k$ is a perfect field and $H$ is cocommutative;

(ii)  $H$ is irreducible cocommutative;

(iii)  The weak action of $H$ on $R$ is inner;

 (iv) $H=(kG)^*$, where $G$ is a finite group;

(v) $H$ is cocommutative.

 Part (i) (ii) are  due to W. Chin
 \cite[Theorem 2, Corollary 1]{Ch92}.
Part (iii) is due to B.J. Blattner and S. Montgomery \cite
[Theorem 2.7]{BM89}. Part (iv) is due to M. Cohen and S.
Montgomery \cite[Theorem 2.9] {CM84b}. Part (v) is due to S.
Montgomery and H.J. Schneider \cite [Corollary 7.13]{MS95}.

         If $H= (kG)^*$, then  relation (\ref {e (2)}) holds, due to
         M. Cohen and
 S. Montgomery \cite[Theorem 4.1] {CM84b}

In  this paper,   we obtain      the relation between $H$-radical
of $H$-module algebra $R$ and radical of $R \# H$. We give
 some sufficient conditions for  (\ref {e (2)})
and (\ref {e (3)}) and  the formulae, which are similar to (\ref
{e (1)}),
 (\ref {e (2)})  and (\ref {e (3)}) for $H$-prime radical respectively.
 We show that (\ref {e (1)}) holds for any
 Hopf algebra $H.$
 Using radical theory and the conclusions in \cite {MS95},
 we also obtain
that    if $H$ is a finite-dimensional semisimple, cosemisimle and
either commutative
   or cocommutative Hopf algebra, then $R$ is $H$-semiprime iff
      $R$ is semiprime iff
$R\#_{\sigma}H$ is semiprime.

In this paper, unless otherwise stated, let $k$ be a field, $R$ be
an algebra with unit over $k$,
 $H$ be a Hopf algebra over $k$ and $H^*$ denote the dual space of $H$.

 $R$ is called a twisted $H$-module algebra if
 the following
conditions are satisfied:

   (i) $H$ weakly acts on  $R$;

(ii) $R$ is a twisted $H$-module, that is,
 there exists a linear map
 $\sigma \in Hom_k (H \otimes H, R) $ such that
$h \cdot (k \cdot r) = \sum \sigma (h_1,k_1)(h_2k_2 \cdot r)
\sigma ^{-1}(h_3,k_3)$ for all $h, k \in H$ and $r \in R$.

  It is clear that if $\sigma $  is trivial, then
   twisted $H$-module algebra $R$ is an $H$-module algebra.

Set
$$ Spec(R) = \{ I \mid I \hbox { \ is a prime ideal of \ } R \};   $$
 $$H\hbox {-}Spec(R) =
\{ I \mid I \hbox { \ is an \ } H\hbox{-prime ideal of } \ R \}.$$

 \section {The Baer radical of twisted  H-module algebras }

In this section, let $k$ be a commutative associative ring with
unit, $H$ be an algebra
 with unit and comultiplication $\bigtriangleup$,
  $R$ be an algebra
over $k$ ($R$ may be without unit) and
 $R$ be  a twisted  $H$-module algebra.

\begin {Definition} \label {10.1.1}
$r_{Hb}(R) := \cap \{I \mid I$ is  an $H$-semiprime ideal of $R$
\};

$r_{bH}(R):= (r_b(R):H)$

$r_{Hb}(R)$ is called the $H$-Baer radical ( or $H$-prime radical
) of twisted $H$-module algebra $R$.

\end {Definition}

\begin {Lemma} \label {10.1.3}
(1) If $E$ is a non-empty subset of $R$,
  then  $(E) = (H \cdot E) + R(H \cdot E) + (H \cdot E)R + R(H \cdot E)R$,
   where $(E)$ denotes the $H$-ideal generated by $E$ in  $R$;

(2) If $I$ is a nilpotent $H$-ideal of $R$, then $I \subseteq
r_{Hb}(R)$.
\end {Lemma}

  {\bf Proof}.
(1) It is trivial.

(2)  If $I$ is a nilpotent $H$-ideal and $P$ is an $H$-semiprime
ideal,
                      then  $(I +P)/P$ is nilpotent simply because
$ (I +P)/P \cong I/(I \cap P) $  \ \ \ \ (as \ algebras) . Thus
$I \subseteq P$ and  $I \subseteq r_{Hb}(R)$. $\Box$

\begin {Proposition}\label {10.1.4}

(1) $r_{Hb}(R)=0$ iff $R$ is $H$-semiprime;

(2) $r_{Hb}(R/r_{Hb}(R))=0$;

(3)  $R$ is $H$-semiprime iff $(H\cdot a)R(H\cdot a) = 0$ always
implies
    $a = 0$ for any $a\in R$;

  $R$ is $H$-prime iff $(H\cdot a)R(H\cdot b) = 0$ always implies
 $a = 0$  or $b = 0$ for any $a$, $b \in R$;

(4) If $R$ is $H$-semiprime, then $W_H(R)=0$.
\end {Proposition}

  {\bf Proof}.
(1) If $r_{Hb}(R)=0$, then $R$ is $H$-semiprime by Lemma \ref
{10.1.3} (2). Conversely, if $R$ is $H$-semiprime, then $0$ is an
$H$-semiprime ideal and so $r_{Hb}(R)=0$ by Definition \ref
{10.1.1}.

(2) If $B/r_{Hb}(R) $ is a nilpotent $H$-ideal of $R/r_{Hb}(R)$,
then $B^k \subseteq r_{Hb}(R)$  for some natural number $k$ and so
$B \subseteq r_{Hb}(R)$, which implies that $R/r_{Hb}(R)$  is
$H$-semiprime. Thus $r_{Hb}(R/r_{Hb}(R))=0$ by part (1).

(3) If $R$ is  $H$-prime and $(H\cdot a)R(H\cdot b) = 0$ for $a$
and $b \in R$, then $(a)^2 (b)^2 = 0$ by Lemma \ref {10.1.3} (1),
where $(a)$ and $(b)$ are the $H$-ideals generated by  $a$ and $b$
in $R$ respectively. Since $R$ is $H$-prime, $a = 0$ or $b= 0$.
Conversely, if both $B$ and $C$ are  $H$-ideals of $R$ and $ BC =
0$, then $(H \cdot a)R(H\cdot b) = 0$  and $a = 0$ or $b = 0$
 for any $ a \in B$ and $b \in C$, which implies  that $B = 0 $ or $C = 0$.
 Thus $R$ is an $H$-prime.
Similarly, the other assertion holds.

(4)    For any $0 \not= a \in R$, there exist
 $b_1 \in R$ and  $h_1, h_1' \in  H$  such that
 $ 0 \not= a_2 = (h_1 \cdot a_1) b_1(h_1' \cdot a_1)  \in
 (H \cdot a_1)R(H \cdot a_1)$ by part (3), where $a_1 = a$.
 Similarly, for $0 \not= a_2 \in R$, there exist
 $b_2 \in A$, $h_2$, $h_2' \in  H$  such that
 $ 0 \not= a_3 = (h_2 \cdot a_2) b_2(h_2' \cdot a_2)  \in
 (H \cdot a_2)R(H \cdot a_2)$, which implies that
 there exists an $H$-$m$-sequence $\{ a_n \}$  such that
  $a_n \not= 0$ for any natural number $n$.
   Thus $W_{H}(R)=0$.  $\Box$

\begin {Theorem} \label {10.1.5}
$r_{Hb}(R)= W_H(R) = \cap  \{ I \mid I { \ is \ an \ } H{
\hbox{-}prime \ ideal \ of \  } R    \} $.
\end {Theorem}

{\bf Proof.} Let  $D = \cap \{ I \mid I $ is an $H$-prime ideal of
$R$ \}. Obviously, $r_{Hb}(R) \subseteq D$.

    If $0 \not= a \not\in W_H(R)$, then there exists
    an $m$-sequence $\{ a_i \}$  in $R$ with
    $a_1=a$ and $a_{n+1} = (h_n \cdot a_n )b_n (h'_n \cdot a_n) \not= 0$
    for $n = 1, 2, \cdots $.
Let  ${\cal F} = \{ I
  \mid I$ is an $H$-ideal of $R$ and $I \cap \{ a_1, a_2, \cdots \}
  = \emptyset \}$.
 By Zorn's Lemma, there exists a
maximal $P$  in ${\cal F}$. If both  $I$ and $J$ are $H$-ideals of
$R$ with  $I \not\subseteq P$ and
 $J \not\subseteq P$
 such that $IJ \subseteq P$, then there exist natural numbers
 $n$  and $m$ such that  $a_n \in I +P$  and
 $a_m \in J+P$. Since $ a_{n +m + 1}
 = (h_{n+m}\cdot a_{n+m})b_{n+m}(h'_{n+m}\cdot a_{n+m}) \in (I+P)(J+P)
 \subseteq P$, we get a contradiction. Thus
 $P$ is an $H$-prime ideal of $R$. Obviously, $a \not\in P$, which implies
 that
  $a \not\in D$.
Therefore $D \subseteq W_H(R)$.

    For any $x \in W_H(R)$,
let $\bar R = R/r_{Hb}(R)$. It follows from
 Proposition \ref {10.1.4} (1) (2) (4) that   $W_H(\bar R)=0$.
 For  any  $H$-$m$-sequence $\{ \bar a_n \}$
 with $\bar a_1 = \bar x$ in  $\bar R$,  there exist
 $\overline b_n \in \overline R$ and $h_n, h_n' \in H$ such that
 $ \overline a_{n+1} =
 (h_n \cdot \overline a_n) \overline b_n (h_n' \cdot \overline a_n)$
                           {~}{~}  for any natural number $n$.
Let   $a_1'=x $  and
 $  a_{n+1}' = (h_n \cdot  a_n')  b_n (h_n' \cdot  a_n')$ {~}{~}
                             for any natural number $n$.
 Since $ \{ a_n'\}$ is an $H$-$m$-sequence with $a_1' = x$ in $R$,
 there exists a natural number $k$ such that $a_k' =0$. It is
 clear that $\overline a_n= \overline {a_n'}$
  for any natural number $n$ by induction. Thus $ \bar a_k = 0$ and
  $\overline x \in W_H(\overline R)$.
  Considering $W_H(\overline R)=0$, we have $x \in r_{Hb}(R),$
  which implies that
  $W_H(R) \subseteq r_{Hb}(R)$. Therefore
  $W_H(R) = r_{Hb}(R)=D$. $\Box$

  \section{ The Baer and Jacobson radicals of crossed products }

By \cite [Lemma 7.1.2] {Mo93}, if $R\#_\sigma H$  is crossed
product defined in \cite [Definition 7.1.1] {Mo93}, then $R$ is a
twisted $H$-module algebra.

Let $R$ be an algebra and $M_{m \times n }(R)$ be the algebra of
$m \times n$ matrices with entries in $R$. For  $ i = 1, 2, \cdots
m $ and  $ j = 1, 2, \cdots n $.  Let $(e_{ij})_{m \times n}$
denote the matrix in $M_{m \times n}(R),$ where $(i,j)$-entry is
$1_R$  and the others are zero. Set

 ${ \cal I}(R) = \{ I \mid I \hbox { is an ideal of } R \}$.

$r_{Hj}(R) := r_j(R\#_{\sigma} H) \cap R$;

$r_{jH}(R):=(r_{j}(R):H)$.

\begin {Lemma}\label{10.2.2}
Let $M_R$ be a free $R$-module with finite rank  and $R'= End
(M_R)$. Then there exists a  unique bijective  map
 $$\Phi: {\cal I}(R) \longrightarrow {\cal I}(R')$$
such that \ \ \  $\Phi (I)M=MI$ and

 (1) $\Phi $ is a map  preserving containments, finite products,
and infinite intersections;

(2) $ \Phi (I) \cong   M_{n \times n}(I)$ for any ideal $I$ of
$R$;

(3) $I$ is a (semi)prime ideal of $R$ iff $\Phi(I)$ is
 a (semi)prime ideal of $R'$;

(4) $r_b(R')= \Phi(r_b(R))$;

(5) $r_j(R')= \Phi(r_j(R))$.
\end {Lemma}

{\bf Proof} Since $M_R$ is a free $R$-module with rank $n$, we can
assume $M=M_{n \times 1}(R)$. Thus $R' = End (M_R)= M_{n \times
n}(R)$ and the module operation of $M$ over $R$ becomes the matrix
operation. Set $M'= M_{1 \times n}(R)$. Obviously, $M'M=R$. Since
$(e_{i1})_{n \times 1} (e_{1j})_{1 \times n} = (e_{ij})_{n \times
n}$  for $i, j = 1, 2, \cdots n$,  $MM' = M_{n \times n}(R) =R'$.
Define
$$\Phi (I) = MIM'$$
\noindent for any ideal $I$ of $R$. By simple computation, we have
that $\Phi (I)$ is an ideal of $R'$ and $\Phi (I)M = MI$. If $J$
is an ideal of $R'$ such that $JM=MI$, then $JM= \Phi (I)M$ and
$J= \Phi(I)$, which implies $\Phi $ is unique. In order to show
that $\Phi$ is a bijection from ${ \cal I}(R)$ onto ${ \cal
I}(R')$, we define  a map  $\Psi$ from ${\cal I}(R')$ to ${ \cal
I}(R)$ sending $I'$ to $M'I'M$ for any ideal $I'$ of $R'$.  Since
$ \Phi \Psi (I')=   MM'I'MM' = I'$ and $\Psi \Phi (I)= M'MIM'M =
I$ for any ideal $I'$ of $R'$ and ideal $I$ of
 $R$,  we have that $\Phi$ is bijective.

(1) Obviously $\Phi$ preserves containments. We see that
 $$\Phi (IJ) = MIJM'= (MIM')(MJM') = \Phi (I) \Phi(J)$$
  for any ideals $I$
and $J$ of $R$. Thus $\Phi$ preserves finite products.
 To show that $\Phi$ preserves infinite
intersections, we first show that
\begin {eqnarray} M (\cap \{ I_\alpha \mid \alpha \in \Omega \}) =
 \cap \{ MI_\alpha \mid \alpha \in \Omega \}  \label {e1.2 (1)}
 \end {eqnarray}
for any  $\{I_\alpha \mid \alpha \in \Omega \} \subseteq {\cal
I}(R)$.
 Obviously, the right side of  relation (\ref {e1.2 (1)})   contains the left side of
  relation (\ref {e1.2 (1)}).
 Let \{$u_1, u_2, \cdots, u_n $ \}  be a basis of $M$ over $R$.
 For any $x \in  \cap \{ MI_{\alpha} \mid \alpha \in \Omega \}$, any
 $\alpha, \alpha' \in \Omega $,
  there exist
 $r_i \in I_\alpha$ and $r_i' \in I_{\alpha'}$ such that
 $x = \sum u_ir_i= \sum u_ir_i'$.
 Since $\{u_i\}$ is a basis, $r_i = r_i'$,  which implies  $x \in
 M (\cap \{ I_\alpha \mid \alpha \in \Omega \})$.
 Thus the  relation (\ref {e1.2 (1)})
 holds.  It follows from  relation (\ref {e1.2 (1)})  that
 $$\Phi(\cap \{ I_\alpha \mid \alpha \in \Omega \})M =
 \cap \{ \Phi(I_\alpha)M \mid \alpha
 \in \Omega \}$$
Since $  \Phi (\cap  \{ I_\alpha  \mid \alpha \in \Omega \})M =
  \cap  \{ \Phi (I_\alpha )M \mid \alpha \in \Omega \}
\supseteq
 (\cap \{ \Phi(I_\alpha) \mid \alpha
 \in \Omega \})M$, we have that
   $$\Phi(\cap \{ I_\alpha \mid \alpha \in \Omega \})
 \supseteq \cap \{ \Phi(I_\alpha) \mid \alpha \in \Omega \}.$$
Obviously,  $$\Phi(\cap \{ I_\alpha \mid \alpha \in \Omega \})
 \subseteq \cap \{ \Phi(I_\alpha) \mid \alpha \in \Omega \}.$$
 Thus           $$\Phi(\cap \{ I_\alpha \mid \alpha \in \Omega \})
 = \cap \{ \Phi(I_\alpha) \mid \alpha \in \Omega \}.$$

(2) Obviously, $ \Phi (I) = MIM' = M_{n \times 1}(R) I M_{1 \times
n}(R) \subseteq M_{n \times n}(I)$. Since $ a(e_{ij})_{n \times n}
= (e_{i 1})_{n \times 1}a (e_{1j})_{1 \times n} \in MIM'$ for all
$a \in I$ and $i, j = 1, 2, \cdots n$,
$$ \Phi (I) = MIM' =
M_{n \times 1}(R) I M_{1 \times n}(R) \supseteq M_{n \times
n}(I).$$ Thus part (2) holds.

(3) Since bijection $\Phi$ preserves  products,  part (3) holds.

(4) We see that
\begin {eqnarray*} \Phi (r_b(R)) &=&
\Phi (\cap \{ I \mid I \hbox { is a prime ideal of }  R \} )  \\
&=& \cap \{ \Phi (I) \mid I \hbox { is a prime ideal of }  R \}
\hbox { by part (1) } \\
&=& \cap \{ \Phi (I) \mid \Phi (I) \hbox { is a prime ideal of }
R' \}
\hbox { by part (3) }   \\
&=& \cap \{ I' \mid I' \hbox { is a prime ideal of }  R' \}
\hbox { since  }  \Phi \hbox { is surjective } \\
&=& r_b(R')
\end {eqnarray*}

(5) We see that
\begin {eqnarray*} \Phi (r_j(R)) &=&
M_{n \times n}(r_j(R)) \hbox { by part (2) }     \\
&=& r_j(M_{n \times n}(R)) \hbox { by \cite [Theorem 30.1] {Sz82} } \\
&=& r_j(R')  { \ . \ \ \ \ \ } \Box
\end {eqnarray*}

     Let $H$ be a finite-dimensional   Hopf algebra and $A = R\#_{\sigma} H$.
Then $A$ is a free right $R$-module  with finite rank by \cite
[Proposition 7.2.11]{Mo93} and $End (A_R) \cong (R \#_{\sigma}
H)\#H^*$ by \cite [Corollary 9.4.17]{Mo93}. By part (a) in the
proof of \cite [Theorem 7.2] {MS95}, it follows that $\Phi$ in
Lemma 1.2 is the same as  in \cite [Theorem 7.2] {MS95}.

 \begin {Lemma} \label {10.2.3} Let $H$ be a finite-dimensional
      Hopf algebra and $A = R\#_{\sigma} H$. Then

(1) If $P$ is an $H^*$-ideal of $A$, then $P=(P\cap R) \#_\sigma
H$.

(2) $ \Phi (I) = (I \#_{\sigma}H)\#H^*$ for  every  $H$-ideal $I$
of $R$;

(3)         \begin {eqnarray}
   \{ P \mid P \hbox { is an } H \hbox {-ideal of } A \#H^* \}
&=& \{ (I \#_{\sigma}H)\#H^* \mid I \hbox { is an \ } H\hbox
{-ideal of } R \}
 \label {e1.3 (1)}
 \end {eqnarray}
 \begin {eqnarray}
 \{ P \mid P \hbox { is an } H^*\hbox {-ideal of } A  \}
 &=& \{ I \#_{\sigma }H \mid I \hbox { is an } H\hbox {-ideal of } R \}
\label {e1.3 (2)}
\end {eqnarray}
$ \{ P \mid P \hbox { is an } H\hbox{-prime ideal of } A \#H^* \}$
           \begin {eqnarray}
&=& \{ (I \#_{\sigma}H)\#H^* \mid I \hbox { is an } H\hbox{-prime
ideal of } R \} \label {e1.3 (3)}
\end {eqnarray}

(4)    $H$-Spec$(R) = \{ (I:H) \mid I  \in$ Spec $(R) \}$;

 (5) \begin {eqnarray}    (\cap  \{ I_{\alpha } \mid  \alpha \in \Omega \} : H)
 &=&
\cap \{ (I_{\alpha } : H) \mid \alpha \in \Omega \} \label {e1.3
(4)}
\end {eqnarray}
 where $I_{\alpha}$  is an ideal
of $R$ for all $\alpha \in \Omega ; $

 (6)  \begin {eqnarray}   ( \cap  \{ I_{\alpha } \mid \alpha \in
 \Omega \} )\#_{\alpha}H &=&
\cap \{ (I_{\alpha } \#_{\sigma } H) \mid \alpha \in \Omega \}
\label {e1.3 (5)}
 \end {eqnarray}
 where $I_{\alpha}$  is an $H$-ideal
of $R$ for all $\alpha \in \Omega $;

(7) $\Phi(r_b(R))= r_b(A\#H^*)$;

(8)  $\Phi(r_j(R))= r_j(A\#H^*)$;

(9) $\Phi(r_{Hb}(R))= r_{Hb}(A\#H^*) = (r_{Hb}(R) \#_\sigma H) \#
H^*$.
\end {Lemma}

{\bf Proof} (1) By  \cite [Corollary 8.3.11] {Mo93}, we have that
 $P = (P \cap R)A = (P\cap R)\#_\sigma H$.

(2) By the part (b) in the proof of \cite [Theorem 7.2]{MS95},
it follows
 that
 $$\Phi(I)= (I \#_{\sigma}H)\#H^*$$
 for every $H$-ideal $I$ of $R$.

(3)  Obviously, the left side of  relation (\ref {e1.3 (1)})
contains
 the right side of  relation (\ref {e1.3 (1)}). If $P$ is an $H$-ideal of
 $A\#H^*$, then
 $P= (P \cap A) \# H^* = (((P \cap A) \cap R) \#_\sigma H) \# H^*$ by
 part (1), which implies that the right side of  relation (\ref {e1.3 (1)}) contains
 the left side. Thus  relation (\ref {e1.3 (1)}) holds.
 Similarly,  relation (\ref {e1.3 (2)}) holds.
 Now, we show that  relation (\ref {e1.3 (3)}) holds. If $P$
 is an $H$-prime ideal of $A \# H^*$, there exists an $H$-ideal $I$ of $R$
 such that $P = (I \#_\sigma H) \#H^*$ by  relation (\ref {e1.3 (1)}). For any
 $H$-ideals $J$ and $J'$ of $R$ with $JJ' \subseteq I$,
 since $\Phi (JJ') = \Phi (J) \Phi(J') \subseteq \Phi (I) =P$  by Lemma \ref {10.2.2} (1), we have that $\Phi (J) \subseteq \Phi(I)$
or $\Phi(J')\subseteq \Phi(I)$, which implies that $J \subseteq I$
or $J' \subseteq I$ by Lemma \ref {10.2.2}. Thus $I$ is an
$H$-prime ideal of $R$. Conversely, if $I$ is an $H$-prime ideal
of $R$ and $P= (I \#_\sigma H)\#H^*$, we claim that $P$ is an

$H$-prime of $A \#H^*$. For any $H$-ideals $Q$ and $Q'$ of $A
\#H^*$ with  $QQ' \subseteq P$, there exist two $H$-ideals $J$ and
$J'$ of $R$ such that $(J\#_\sigma H)\#H^* = Q$ and $(J'\#_\sigma
H)\#H^* = Q'$ by
 relation (\ref {e1.3 (1)}).
Since $\Phi(JJ') = \Phi(J) \Phi(J') = QQ'\subseteq P = \Phi(I) $,
$JJ' \subseteq I$, which implies $J \subseteq I$, or $J' \subseteq
I$, and so $Q \subseteq P$ or $Q' \subseteq P$. Thus $P$ is an
$H$-prime ideal of $A\#H^*$. Consequently,  relation (\ref {e1.3
(3)}) holds.

 (4) It follows from \cite [Lemma 7.3 (1) (2)] {MS95}.

(5) Obviously,  the right side  of relation (\ref {e1.3 (4)})
contains the left side. Conversely, if $x \in \cap \{(I_\alpha:H)
\mid \alpha \in \Omega \}$, then $x \in (I_{\sigma} :H)$ and
$h\cdot x \in I_\alpha $ for all $\alpha \in \Omega, h \in H$,
which implies that $h \cdot x \in \cap \{ I_\alpha \mid \alpha \in
\Omega \} $ and $x \in (\cap \{I_\alpha \mid \alpha \in \Omega
\}:H)$. Thus  relation (\ref {e1.3 (4)}) holds.

(6) Let $\{ h^{(1)}, \cdots, h^{(n)} \}$  be a basis of $H$.
Obviously,
 the right side of
 relation (\ref {e1.3 (5)}) contains
 the left side of  relation (\ref {e1.3 (5)}). Conversely,
 for   $ u \in \cap \{ (I_{\alpha } \#_{\sigma } H)
 \mid \alpha \in \Omega \}$ and   $\alpha, \alpha' \in \Omega, $
there exist   $r_{i} \in I_\alpha$ and $r_i' \in I_{\alpha'}$
 such that $u = \sum_{i=1}^{n} r_i \# h^{(i)} =
  \sum_{i=1}^{n} r_i' \# h^{(i)}$.
   Since  $\{ h^{(1)}, \cdots, h^{(n)} \}$  is linearly independent,
   we have that $r_i= r_i'$, which implies that
 $ u \in (\cap  \{ I_{\alpha } \mid \alpha \in \Omega \} ) \#_{\alpha} H$.
Thus  relation (\ref {e1.3 (5)}) holds.

(7) and (8) follow from Lemma \ref{10.2.2}(4)(5).

(9)  We see  that
\begin {eqnarray*}
r_{Hb}(A \#H^*) &=& \cap \{ P \mid P \hbox { is an } H\hbox
{-prime ideal of }
 A \#H^*      \} \hbox { \ \  by Theorem \ref {10.1.5} }\\
 &=& \cap \{ (I \#_\sigma H)\#H^*  \mid I \hbox { is an }
 H\hbox {-prime ideal of }  R       \}  \hbox { \ \ \
 by relation (\ref {e1.3 (3)}) } \\
  &=& ( \cap \{ I \#_\sigma H \mid I \hbox { is an } H\hbox {-prime ideal of }
   R   \} ) \# H^*    \hbox { \ \ \ by  part (6) }  \\
 &=& ((\cap \{ I \mid I \hbox { is an }
 H\hbox {-prime ideal of } R \}) \#_\sigma H) \# H^*
  \hbox { \ \ \ by part (6) }  \\
 &=& (r_{Hb}(R) \#_\sigma H) \# H^*  \hbox { \ \ by Theorem \ref {10.1.5} }\\
&=& \Phi (r_{Hb}(R))    \hbox { \ \ \ by part (2) }.
\end {eqnarray*}

  (10) If $H$ is cosemisimple, then $H$ is semisimple  by \cite [Theorem 2.5.2]
{Mo93}. Conversely, if $H$ is semisimple, then $H^*$ is
cosemisimple. By \cite [Theorem 2.5.2]{Mo93}, $H^*$ is semisimple.
Thus $H$ is cosemisimple. $\Box$

\begin {Proposition}\label {10.2.4}
  (1)  $r_{Hb}(R) \subseteq  r_{b}(R\#_{\sigma}H) \cap R \subseteq r_{bH}(R)$;

  (2)  $ r_{Hb}(R)\#_{\sigma}H \subseteq r_{b}(R \#_{\sigma} H)$.
  \end {Proposition}

  {\bf Proof}.
   (1) If $P$ is a prime ideal of $R \#_{\sigma} H$, then $P \cap R$ is
    an $H$-prime    ideal of $R$ by \cite [Lemma 1.6]{Ch91}.
   Thus  $r_b(R \#_{\sigma} H) \cap R= \cap \{ P\cap R \mid
   P$ is prime ideal of $R \#_{\sigma} H \}  \supseteq r_{Hb}(R)$.
For any $a \in r_b(R \#_\sigma H) \cap R$ and any $m$-sequence $
\{  a_i \} $ in $R$ with $a_1 = a $, it is easy to check that $\{
a_i \}$ is also an $m$-sequence in $R \#_\sigma H$. Thus $a_n = 0
$  for some natural $n$, which implies $a \in r_b(R)$. Thus
$r_b(R\#_\sigma H) \cap R \subseteq r_{bH}(R)$ by \cite [Lemma
1.6]{BM89}

(2) We see that
\begin {eqnarray*}
    r_{Hb}(R) \#_{\sigma} H &=&    (r_{Hb}(R) \#_{\sigma} 1)(1 \#_{\sigma} H)\\
    &\subseteq &  r_b(R \#_{\sigma} H) ( 1 \#_\sigma H)
    \hbox { \ \  by  part  (1)}\\
    &\subseteq & r_b(R \#_{\sigma} H). { \ \ \ \ }  \Box
    \end {eqnarray*}

 \begin {Proposition} \label {10.2.5} Let $H$ be finite-dimensional
      Hopf algebra and $A = R\#_{\sigma} H$.  Then

(1) $r_{H^*b}(R\#_{\sigma}H) =  r_{Hb}(R)\#_{\sigma}H$;

(2) $r_{Hb}(R)= r_{bH}(R)= r_b(R\#_{\sigma}H) \cap R$.
\end {Proposition}

{\bf Proof}

(1) We see that
\begin {eqnarray*}
r_{H^*b}(R \#_\sigma H) &=& \cap \{ P \mid P
\hbox { is an } H^*\hbox {-prime ideal of } A \} \\
&=& \cap \{ I \#_\sigma H \mid I \hbox { is an } H \hbox {-prime
ideal of } R \}  \hbox
{ \ \ \ ( by  \cite [Lemma 7.3 (4)] {MS95} ) }\\
&=& ( \cap \{ I \mid I \hbox { is an } H \hbox {-prime ideal of }
R \}) \#_\sigma H
\hbox { \ \ \ ( by Lemma \ref {10.2.3} (6)) }\\
&=& r_{Hb}(R) \#_\sigma H.
\end {eqnarray*}

 (2) We see that
\begin {eqnarray*} r_{Hb}(R) &=& \cap \{ P \mid P { \ is \ an \ }
H \hbox {-prime  ideal  of  } R \} \\
   &=& \cap \{(I:H) \mid I \in { \ Spec} (R) \}
   \hbox { \ by  Lemma \ref {10.2.3} \  part (4) }\\
   &=& ( \cap \{ I \mid I \in {Spec}(R) \}:H  )
 \hbox   { \ by \ Lemma \ref {10.2.3} \  part (5) }  \\
   &=& (r_b(R):H) \\
   &=& r_{bH}(R).
   \end{eqnarray*}

 Thus it follows from Proposition \ref {10.2.4}(1)
   that $r_{Hb}(R)= r_b(R\#_{\sigma}H)\cap R = r_{bH}(R)$.
$\Box$

      \begin {Theorem} \label {10.2.6}.
Let $H$ be a finite-dimensional Hopf algebra and  the weak  action
of $H$ be  inner. Then

(1) $r_{Hb}(R)= r_b(R)= r_{bH}(R)$;

 Moreover, if $H$ is semisimple, then

 (2) $r_{b}(R\#_{\sigma} H)= r_{Hb}(R)\#_{\sigma} H$.
 \end {Theorem}

{\bf Proof} (1) Since the weak action is inner, every ideal of $R$
is an $H$-ideal, which implies that
 $r_{Hb}(R)= r_b(R) = r_{bH}(R)$ by Proposition \ref {10.2.5} (2).

(2) Considering  Proposition \ref {10.2.4}(2), it suffices to show
 $r_{b}(R\#_{\sigma} H)\subseteq  r_{Hb}(R)\#_{\sigma} H.$
It is clear that
\begin {eqnarray} (R\#_{\sigma} H)/ (r_{Hb}(R)\#_{\sigma} H)
\cong (R/r_{Hb}(R))\#_{\sigma} H  \hbox { \ \ \ \ (  \ as \
algebras )}. \label {e1.6 (1)}
\end {eqnarray}
It follows by \cite [Theorem 7.4.7]{Mo93} that $
(R/r_{Hb}(R))\#_{\sigma} H$
 is semiprime.
Therefore $$r_{b}(R \#_{\sigma} H) \subseteq r_{Hb}(R)\#_{\sigma}
H. \ \ \Box$$

   \begin {Theorem}\label {10.2.7} Let $H$ be  a finite-dimensional, semisimple
   and either commutative or cocommutative Hopf
algebra and let $A = R\#_{\sigma} H$. Then

(1) ${r_{b}(R\#_{\sigma}H) =  r_{Hb}(R)\#_{\sigma}H};$

(2) $R$ is $H$ semiprime iff $R \#_\sigma H$ is semiprime.

Moreover, if $H$ is cosemisimple, or char $k$ does not divide dim
$H$, then  both part (3) and part (4) hold:

(3) $r_{Hb}(R)= r_{bH}(R)= r_b(R)$;

(4) $R$ is $H$-semiprime iff $R$ is semiprime iff $R\#_{\sigma}H$
is semiprime.

\end {Theorem}
{\bf Proof} (1) Considering Proposition \ref {10.2.4}(2), it
suffices to show
 $$r_{b}(R\#_{\sigma} H)\subseteq  r_{Hb}(R)\#_{\sigma} H.$$
It follows by \cite [Theorem 7.12 (3)]{MS95} that $
(R/r_{Hb}(R))\#_{\sigma} H$
 is semiprime. Using relation (\ref {e1.6 (1)}), we have
 that $r_{b}(R \#_{\sigma} H) \subseteq r_{Hb}(R)\#_{\sigma} H$.

(2)  It follows from part (1) and Proposition \ref {10.1.4} (1).

(3)  By \cite [Theorem 4.3 (1)] {La71}, we have that $H$ is
semisimple and cosemisimple.

 We see  that
 \begin{eqnarray*}
\Phi(r_b(R)) &=& r_{b}(A \# H^*) \hbox { \ \ \ by \ Lemma \ref {10.2.3} \ (7)}\\
 &=& r_{H^*b}(A) \# H^*    \hbox { \ \ \  by \  part \  (1) } \\
 &=& (r_{Hb}(R) \#_{\sigma} H) \# H^* \hbox { \ \ \ by \ Proposition \ref {10.2.5} (1) } \\
 &=& \Phi (r_{Hb}(R))  \hbox { \ \ \ by \ Lemma \ \ref {10.2.3} (2). }
 \end {eqnarray*}
 Thus $r_b(R)= r_{Hb}(R).$

(4) It immediately follows from  part (2) and part (3). $\Box$

We now provide an example to show that the Baer radical $r_b(R)$
 of $R$ is not $H$-stable  when $H$ is not cosemisimple.

Example: Let $k$ be a field of characteristic $p > 0$ and $R =
k[x]/(x^p)$. Then we can define a derivation $d$ on $R$ by sending
$x$  to $x+1$. Then $d^2(x)=d(x+1)=d(x)$ and then, by induction,
$d^p(x)=d(x)$. It follows that $d^p=d$ on all of $R$.Thus
$H=u(kd)$, the restricted enveloping  algebra, is semisimple by
\cite [Theorem 2.3.3]{Mo93}.
 clearly $H$ acts on $R$, but $H$
does not stabilize the Baer radical of $R$ which is the principal
ideal generated by $x$. Note also that $H$ is commutative  and
cocommutative.

\begin {Proposition}\label {10.2.8}
If $R$ is an $H$-module algebra, then
$$r_{Hj}(R) = \cap \{(0:M)_R \mid M { \ is \ an \ irreducible \  }
R\hbox {-}H\hbox {-} { module} \}.$$ That is, $r_{Hj}(R)$ is the
$H$-Jacobson radical of the $H$-module algebra $R$ defined in
\cite {Fi75}.
\end {Proposition}

{\bf Proof.}  It is easy to show that $M$ is an irreducible
$R$-$H$-module iff $M$ is an irreducible $R\#H$-module by \cite
[Lemma 1]{Fi75}. Thus
\begin {eqnarray*}
r_{Hj}(R) &=& r_{j}(R\#H)\cap R \hbox { \ by \ definition \ref {10.1.1}} \\
&=& (\cap \{ (0:M)_{R\#H} \mid M \hbox { \ is \ an \ irreducible \
}
R\# H \hbox {-module} \} ) \cap R \\
&=& \cap \{ (0:M)_R \mid M  \hbox { \ is \ an \ irreducible \ }
R\hbox {-}H \hbox {-module} \}. \Box
\end {eqnarray*}

\begin {Proposition} \label {10.2.9}

  (1)  $ r_{j}(R\#_{\sigma}H ) \cap R =r_{Hj}(R) \subseteq r_{jH}(R); $

  (2)  $r_{Hj}(R)\#_{\sigma} H \subseteq r_{j}(R \#_\sigma H)$.
  \end {Proposition}

  {\bf Proof}.
  (1)  For any $a \in r_j(R \#_{\sigma} H) \cap R$, there exists
$u = \sum_{i} a_i \# h_i \in R \#_{\sigma} H$ such that
$$ a + u + au = 0. $$
Let $(id \otimes \epsilon )$ act on the above equation. We get
that $a + \sum a_i\epsilon (h_i) + a (\sum a_i \epsilon (h_i))=0$,
which implies that $a$ is a right quasi-regular element in $R$.
Thus $r_j(R \#_{\sigma} H) \cap R \subseteq r_{jH}(R)$.

   (2)  It is similar to the proof of Proposition \ref {10.2.4} (2).
   $\Box$

 \begin {Proposition} \label {10.2.10} Let $H$ be a finite-dimensional
      Hopf algebra and $A = R\#_{\sigma} H$. Then

  (1) $r_{jH}(R) \#_{\sigma} H = r_{H^*j}(R \#_\sigma H);$

(2)  $r_{Hj}(R)=r_{jH}(R);$

  (3) $r_{Hj}(A \# H^*) = (r_{Hj}(R) \#_\sigma H) \# H^*$.
  \end {Proposition}

{\bf Proof} (1)         We see that
\begin{eqnarray*}
(r_{jH}(R)\#_{\sigma} H)\#H^* &=& \Phi(r_{jH}(R))\\
 &=& (\Phi(r_j(R) )\cap A) \#H^*  \hbox { \ \  by\   \cite [Theorem 7.2]{MS95} }\\
&=&(r_j(A\#H^*)\cap A) \# H^* \hbox {  \  \ by\  Lemma \ref {10.2.3} (8) } \\
&=& r_{H^*j}(A)\#H^* \hbox { \ by \ Definition \ref{10.1.1} }.
\end {eqnarray*}
Thus $r_{H^*j}(A)= r_{jH}(R)\#_{\sigma} H.$

(2) We see that
\begin {eqnarray*}
 r_{Hj}(R) &=& r_j(A) \cap R \\
& \supseteq & r_{H^*j}(A) \cap R \hbox  { \ by \ Proposition \ref {10.2.9} (1) }\\
&=&  r_{jH}(R)  \hbox { \ by \ part  \  (1) }.
\end {eqnarray*}
It follows by Proposition \ref{10.2.2} (1) that $r_{Hj}(R)=
r_{jH}(R)$.

(3) It immediately follows from part (1) (2). $\Box$

By Proposition \ref {10.2.9} and \ref {10.2.10}, it is clear that
if $H$ is a finite-dimensional Hopf algebra, then
 relation (\ref {e (2)})  holds iff  relation (\ref {e (3)}) holds.

 \begin {Theorem} \label {10.2.11}
 Let $H$ be a finite-dimensional  Hopf algebra and
  the weak action of $H$ be inner.
Then

(1) $r_{Hj}(R)= r_j(R)= r_{jH}(R).$

Moreover, if $H$ is semisimple, then

(2) $r_{j}(R\#_{\sigma} H)= r_{Hj}(R)\#_{\sigma} H$.
\end {Theorem}

{\bf Proof} (1) Since the weak action is inner, every ideal of $R$
is an $H$-ideal and $r_j(R) = r_{jH}(R)$. It follows from
Proposition \ref {10.2.3}(2) that $r_{Hj}(R)= r_{jH}(R)= r_j(R)$.

(2) Considering Proposition \ref {10.2.9}(2),
 it suffices to show
 $$r_{j}(R\#_{\sigma} H)\subseteq  r_{Hj}(R)\#_{\sigma} H.$$
It is clear that  $$(R\#_{\sigma} H)/ (r_{Hj}(R)\#_{\sigma} H)
\cong (R/r_{Hj}(R))\#_{\sigma} H  \hbox { \ \ \ (as algebras). }$$
It follows by \cite [Corollary 7.4.3]{Mo93}  and part (1) that $
(R/r_{Hj}(R))\#_{\sigma} H$  is semiprimitive. Therefore $r_{j}(R
\#_{\sigma} H) \subseteq r_{Hj}(R)\#_{\sigma} H$. $\Box$

\begin {Theorem}\label {10.2.12} Let $H$ be  a finite-dimensional, semisimple Hopf
algebra, let $k$ be an algebraically closed field and let $A =
R\#_{\sigma} H$. Assume  $H$ is cosemisimple or char $k$ does not
divide dim $H$.

 (1) If $H$ is cocommutative, then
 $$r_{Hj}(R)= r_{jH}(R) = r_j(R);$$

(2) If $H$ is commutative, then
$$r_{j}(R \#_{\sigma}H)= r_{Hj}(R) \#_{\sigma}H.$$
\end {Theorem}

{\bf Proof}     By \cite [Theorem 4.3 (1)] {La71} , we have that
$H$ is semisimple and cosemisimple.

(1)  If $g \in G(H)$, then the weak action of $g$ on $R$  is an
algebraic homomorphism, which implies that $g \cdot r_j(R)
\subseteq r_j(R)$. Let $H_0$ be the coradical of $H, H_1 = H_0
\wedge H_0, H_{i+1} = H_0 \wedge H_i$ for $i = 1, \cdots, n$,
where $n$ is the dimension $dimH$ of $H$. It is clear that
$H_0=kG$ with $G =G(H)$ by \cite [Theorem 8.0.1 (c)]{Sw69a} and
$H= \cup H_i$. It is easy to show that if $k > i $, then  $$H_i
\cdot (r_j(R))^k\subseteq r_j(R)$$ by induction for $i$. Thus
$$H \cdot (r_j(R))^{dimH+1} \subseteq r_j(R),$$
which implies that $(r_j(R))^{dimH +1}\subseteq r_{jH}(R)$.

We see that
\begin{eqnarray*}
 r_{j}(R/r_{jH}(R)) &=& r_j(R)/r_{jH}(R)\\
 &=& r_b(R/r_{jH}(R)) \hbox { \ since \ }  r_j(R)/r_{jH}(R) \hbox { \ is
 \ nilpotent } \\
 &=& r_{bH}(R/r_{jH}(R)) \hbox { \ by \ Theorem \ \ref {10.2.7} (3)} \\
 &\subseteq &  r_{jH}(R/r_{jH}(R)) \\
 &=& 0 \ \ .
\end {eqnarray*}

Thus       $r_{j}(R) \subseteq r_{jH}(R)$, which implies that
  $r_j(R)=r_{jH}(R)$.

(2) It immediately follows from part (1) and Proposition
 \ref {10.2.10}(1) (2).
$\Box$

\section {The general theory of $H$-radicals for twisted $H$-module algebras }

In this section we give the general theory of $H$-radicals for
twisted $H$-module algebras.

\begin {Definition}\label {10.3.1}
Let $r$ be a property of $H$-ideals of twisted $H$-module
algebras. An $H$- ideal $I$ of twisted $H$-module algebra $R$ is
called   an $r$-$H$-ideal of $R$ if it is of the $r$-property. A
twisted $H$-module algebra $R$  is called   an $r$-twisted
$H$-module algebra if it is $r$-$H$-ideal of itself. A property
$r$ of $H$-ideals of twisted $H$-module algebras   is called an
 $H$-radical property if the following conditions are satisfied:

(R1) Every twisted $H$-homomorphic image of $r$-twisted $H$-module
algebra is an $r$ twisted $H$-module algebra;

(R2) Every twisted $H$-module algebra $R$ has  the maximal
$r$-$H$-ideal $r(R)$;

(R3)  $R/r(R)$ has not any non-zero $r$-$H$-ideal.

\end {Definition}
We call $r(R)$ the $H$-radical of $R$.

\begin {Proposition}\label {10.3.2}
Let $r$ be an ordinary hereditary radical property for rings. An
$H$-ideal $I$ of twisted $H$-module algebra $R$ is called  an
$r_H$-$H$-ideal  of $R$ if $I$ is an $r$-ideal of ring $R$. Then
$r_H$ is an $H$-radical property for twisted $H$-module algebras.
\end {Proposition}

{\bf Proof.} (R1). If $(R,\sigma )$ is an $r_H$-twisted $H$-module
algebra
 and $(R, \sigma )  \stackrel {f} { \sim }  (R', \sigma ')$, then $r(R') =R'$
 by ring theory. Consequently, $R'$  is an $r_H$-twisted $H$-module algebra.

(R2). For any  twisted $H$-module algebra $R$, $r(R)$  is the
maximal $r$-ideal of $R$ by ring theory. It is clear that $r(R)_H$
is the maximal $r$-$H$-ideal,  which is an $r_H$-$H$-ideal of $R$.
Consequently, $r_H(R) = r(R)_H$ is the maximal $r_H$-$H$-ideal of
$R$.

(R3). If $I/r_H(R)$  is an $r_H$-$H$-ideal of $R/r_H(R)$, then $I$
is an $r$-ideal of algebra $R$ by ring theory. Consequently, $I
\subseteq r(R)$ and $I \subseteq r_H(R).$ $\Box$

\begin {Proposition}\label {10.3.3}
$r_{Hb}$ is an $H$-radical property.
\end {Proposition}

{\bf Proof.} (R1). Let $(R,\sigma )$ be an $r_{Hb}$-twisted
$H$-module algebra
 and $(R, \sigma )  \stackrel {f} { \sim }  (R' \sigma ')$.
    For any $x' \in R'$
and  any  $H$-$m$-sequence $\{  a_n' \}$    in $R'$
 with $ a_1' =  x' $,  there exist
 $ b_n' \in  R'$ and $h_n, h_n' \in H$ such that
 $  a_{n+1}' =
 (h_n \cdot  a_n')  b_n' (h_n' \cdot  a_n')$
                           {~}{~}  for any natural number $n$.
Let   $a_1, b_i \in R $  such that $f(a_1) =x'$ and
$f(b_i)=b_i'$ for $i = 1, 2, \cdots .$  Set
 $  a_{n+1} = (h_n \cdot  a_n)  b_n (h_n' \cdot  a_n)$ {~}{~}
                             for any natural number $n$.
 Since $ \{ a_n \}$ is an $H$-$m$-sequence in $R$,
 there exists a natural number $k$ such that $a_k =0$. It is
 clear that $f( a_n)= a_n'$
  for any natural number $n$ by induction. Thus $  a_k' = 0,$
  which implies that $x' $  is an $H$-$m$-nilpotent element.
Consequently, $R'$ is an $r_{Hb}$-twisted $H$-module algebra.

(R2). By \cite [Theorem 1.5] {Zh98a}, $r_{Hb}(R)= W_H(R) = \{ a
\mid a$ is an $H$-$m$-nilpotent element in $R \}$. Thus
$r_{Hb}(R)$ is the maximal $r_{Hb}$-$H$-ideal of $R$.

(R3).  It immediately follows from  \cite [Proposition 1.4]
{Zh98a}. $\Box$

\section {The relations among  radical of $R$  , radical of
  $R \# _\sigma H$, and $H$-radical of $R$ }

In this section we give the relation among the Jacobson radical
$r_j (R)$ of $R$  ,the Jacobson
 radical $r_j(R \# _\sigma H)$ of
  $R \# _\sigma H$, and $H$-Jacobson radical $r_{Hj}(R)$ of $R$.

In this section, let $k$ be a field, $R$ an algebra with unit, $H$
a Hopf algebra over $k$ and $R \# _\sigma H$ an algebra with unit.
Let $r$ be a hereditary radical property for rings which satisfies
$$r(M_{n \times n} (R)) = M_{n \times n} (r(R))$$ for any twisted
$H$-module algebra $R$.

 Example.  $r_j$, $r_{bm}$ and $r_n$ satisfy the above conditions
by \cite {Sz82}. Using \cite [Lemma 2.1 (2)]{Zh98a},we can easily
prove that $r_b$ and $r_l$   also satisfy the above conditions.

\begin {Definition} \label {10.4.1}
$\bar r_H (R) := r(R \# _\sigma H) \cap R$  and $r_H (R) :=
(r(R):H)$.
\end {Definition}

     If $H$ is a finite-dimensional   Hopf algebra and $M = R\#_{\sigma} H$,
     then $M$ is a free right $R$-module  with finite rank by
\cite [Proposition 7.2.11]{Mo93} and $End (M_R) \cong (R
\#_{\sigma} H)\#H^*$ by \cite [Corollary 9.4.17]{Mo93}. It follows
from  part (a) in the proof of \cite [Theorem 7.2] {MS95} that
 there exists a  unique bijective  map
 $$\Phi: {\cal I}(R) \longrightarrow {\cal I}(R')$$
such that \ \ \  $\Phi (I)M=MI,$ where $R' = (R \#_\sigma H)\#
H^*$  and

 ${ \cal I}(R) = \{ I \mid I \hbox { is an ideal of } R \}$.

\begin {Lemma}\label {10.4.2}  If $H$ is a finite-dimensional Hopf algebra, then

  $ \Phi (r(R)) = r((R \# _\sigma H) \# H^*).$

\end {Lemma}
{\bf Proof.}  It is similar to the proof of \cite [lemma 2.1 (5)]
{Zh98a}. $\Box$

\begin {Proposition}\label {10.4.3}
$\bar r_H (R) \# _\sigma H \subseteq r_{H^*}(R \# _\sigma H)
                           \subseteq r(R \# _\sigma H). $

\end {Proposition}
{\bf Proof.}
 We see that
\begin {eqnarray*}
   \bar r_{H}(R) \#_{\sigma} H &=&    (\bar r_{H}(R) \#_{\sigma} 1)(1 \#_{\sigma} H)\\
    &\subseteq &  r(R \#_{\sigma} H) ( 1 \#_\sigma H) \\
    &\subseteq & r(R \#_{\sigma} H). { \ \ \ \  }
\end {eqnarray*}
Thus      $\bar r_H (R) \# _\sigma H \subseteq r_{H^*}(R \#
_\sigma H)$ since     $\bar r_H(R) \# _\sigma H $  is an
$H^*$-ideal of $R \# _\sigma H $.   $\Box$

\begin {Proposition}\label {10.4.4}
If $H$ is a finite-dimensional
      Hopf algebra,  then

(1) $r_H (R) \# _\sigma H =  \bar r_{H^*}(R \# _\sigma H)$;

Furthermore, if               $\bar r_H \le r_H $, then

(2)    $\bar r_H = r_H$ and $r_H (R) \# _\sigma H \subseteq r(R \#
_\sigma H)$;

(3)    $R\# _\sigma H$  is $r$-semisimple for any $r_H$-semisimple
$R$ iff
   $$r(R \#_\sigma H) = r_H(R) \# _\sigma H.$$

\end {Proposition}
{\bf Proof.} Let $A = R \# _\sigma H.$

(1)         We see that
\begin{eqnarray*}
(r_{H}(R)\#_{\sigma} H)\#H^* &=& \Phi(r_{H}(R))\\
 &=& (\Phi(r(R) )\cap A) \#H^*  \hbox { \ \  by\   \cite [Theorem 7.2]{MS95} }\\
&=&(r(A\#H^*)\cap A) \# H^* \hbox {  \  \ \ by  Lemma \ref {10.4.2}  } \\
&=& \bar r_{H^*}(A)\#H^* \hbox { \ by \ Definition \ref{10.4.1} }.
\end {eqnarray*}
Thus $ \bar r_{H^*}(A)= r_{H}(R)\#_{\sigma} H.$

(2) We see that
\begin {eqnarray*}
\bar  r_{H}(R) &=& r(A) \cap R \\
& \supseteq & \bar r_{H^*}(A) \cap R \hbox  { \ by assumption }\\
&=&  r_{H}(R)  \hbox { \ by \ part  \  (1) }.
\end {eqnarray*}
Thus  $\bar r_{H}(R)= r_{H}(R)$ by assumption.

(3) Sufficiency is obvious. Now we show the necessity.
            Since $$r((R \# _\sigma H)/(r_H (R) \#_\sigma H))
            \cong r(R/r_H(R) \#_{\sigma '}H) =0 ,$$
            we have $r (R \#_\sigma H) \subseteq r_H(R) \#_\sigma H.$
Considering part (2), we have

   $$r(R \#_\sigma H) = r_H(R) \# _\sigma H. { \ \ \ } \Box  $$

\begin {Corollary}\label {10.4.5}
Let $r$ denote $r_b, r_l,  r_j, r_{bm}$ and $r_n.$ Then

(1) $\bar r_H \le r_H;$

Furthermore, if $H$ is a finite-dimensional Hopf algebra, then

(2) $\bar r_H = r_H $;

(3) $R \# _\sigma H $ is $r$- semisimple for any
$r_{H}$-semisimple $R$  iff $r(R \# _\sigma H)= r_{H}(R)\#_\sigma
H$;

(4)  $R \# _\sigma H $ is $r_j$- semisimple for any
$r_{Hj}$-semisimple $R$  iff $r_j(R \# _\sigma H)=
r_{Hj}(R)\#_\sigma H.$

\end {Corollary}
{\bf Proof.}  (1) When $r= r_b$  or $r = r_j$, it has been proved
in \cite [Proposition 2.3 (1) and 3.2 (1)] {Zh98a} and in the
preceding sections. The others can similarly be proved.

(2) It follows from Proposition \ref {10.4.4} (2).

(3) and (4)  follow from part (1) and Proposition \ref {10.4.4}
(3). $\Box$
\begin {Proposition}\label {10.4.6}
 If $H = kG$ or  the weak action of $H$ on $R$ is inner,
then

(1). $r_H (R) = r(R)$;

(2)  If, in  addition,   $H$ is a finite-dimensional
 Hopf algebra and $\bar r _H \le r _H$, then
 $r_H(R) = \bar r_H(R) = r(R).$

\end {Proposition}

{\bf Proof.} (1)  It is trivial.

(2) It immediately follows from part (1) and Proposition \ref
{10.4.1} (1) (2).  $\Box$

\begin {Theorem}\label {10.4.8}
Let $G$ be a finite group and $\mid G \mid ^{-1} \in k$. If $H=
kG$ or $H= (kG)^*$, then

(1)  $r_j (R) = r_{Hj}(R)= r_{jH}(R)$;

  (2)   $r_{j}(R\# _\sigma  H)= r_{Hj}(R)\# _\sigma H.$
                  \end {Theorem}
{\bf Proof.} (1) Let $H =kG$.  We  can easily check $r_j(R) =
r_{jH}(R)$  using the method similar to  the proof of
 \cite [Proposition 4.6]  {Zh97b}. By  \cite [Proposition 3.3 (2) ] {Zh98a},
 $r_{Hj}(R) =r_{jH}(R)$. Now,
we only need to show that
  $$r_j (R) = r_{H^*j}(R). $$
 We see that
\begin {eqnarray*}
  r_{j}((R\# _\sigma  H^*) \# H) &=&
   r_{H^*j}((R\# _\sigma H^*) \# H)   \hbox { \ \ by \cite [Theorem 4.4 (3)]
   {CM84b} } \\
   &=&  r_{Hj}(R\# _\sigma  H^*) \# H  \hbox { \ \ by \cite [ Proposition
   3.3 (1)] {Zh98a}}  \\
   &=& (r_{H^*j}(R)\# _\sigma H^* ) \# H  \hbox { \ \ by \cite [ Proposition
   3.3 (1)] {Zh98a}}.
\end {eqnarray*}
On the one hand, by       \cite [ Lemma 2.2 (8)] {Zh98a}, $\Phi
(r_j (R)) =  r_{j}((R\# _\sigma  H^*) \# H).$ On the other hand,
we have that $ \Phi (r_{H^*j} (R)) =  ( r_{H^*j}(R )\# _\sigma
H^*) \# H $  \ \
 by \cite [Lemma 2.2 (2)]    {Zh98a}.
 Consequently, $r_j (R) = r_{H^*j}(R)$.

(2) It immediately follows from part (1) and \cite [Proposition
3.3 (1) (2)] {Zh98a}.
  $\Box$

\begin {Corollary}\label {10.4.9}
Let $H$ be a semisimple and cosemisimple Hopf algebra over
algebraically closed field $k$. If $H$ is commutative or
cocommutative, then
  $$r_j (R) = r_{Hj}(R)= r_{jH}(R)  \hbox { \ \ \ \ and  \ \ \ }
  r_{j}(R\# _\sigma  H)= r_{Hj}(R)\# _\sigma H.$$
                  \end {Corollary}
{\bf Proof.} It immediately follows from Theorem \ref {10.4.8} and
\cite [Lemma 8.0.1 (c)] {Sw69a}. $\Box$

    We give an example to show that conditions in Corollary \ref {10.4.9} can not
be omitted.

\begin {Example}\label {10.4.10} (see \cite [Example P20] {Zh98a})
Let $k$ be a field of characteristic $p > 0 $, $R = k[x]/(x^p)$.
We can define a derivation on $R$ by sending $x$ to $x +1$. Set
$H=u(kd)$, the restricted enveloping algebra, and $A = R\#H.$ Then

(1) $r_b (A \# H^*) \not= r_{H^*b}(A) \# H^*$;

(2) $r_j (A \# H^*) \not= r_{H^*j}(A) \# H^*$;

(3) $r_j (A \# H^*) \not\subseteq  r_{jH^*}(A) \# H^*$.

\end {Example}
{\bf Proof.}        (1) By \cite [Example P20]  {Zh98a}, we have
$r_b(R) \not=0$ and $r_{bH}(R)=0.$ Since $\Phi (r_b(R) ) = r_b(A
\#H^*)\not=0$  and $\Phi (r_{bH}(R))= r_{bH^*}(A) \#H^*=0,$ we
have that  part (1) holds.

(3) We see that $r_j(A \# H^*) = \Phi (r_j(R))$  and $r_{Hj}(A)
\#H^*= \Phi (r_{Hj}(R)).$  Since $R$ is commutative, $r_j(R)=
r_b(R).$  Thus $r_{Hj}(R)= r_{jH}(R)= r_{bH}(R)=0$ and $r_j(R) =
r_b(R) \not=0$, which implies
 $r_j(A\#H^*) \not\subseteq r_{jH^*}(A)\#H^*$.

(2) It follows from part (3).
 $\Box$

This example also answer the question J.R. Fisher asked in \cite
{Fi75} :
$$\hbox {Is  \ \ \ } r_j (R \# H) \subseteq  r_{jH}(R) \# H \hbox { \ \ \  ?} $$

If $F$ is an extension field of $k$, we write  $R^F$  for $R
\otimes _k F$ (see \cite [P49 ]{MS95}) .

\begin {Lemma}\label {10.4.11}  If $F$ is an extension field  of $k$, then

(1)  $H$ is a semisimple  Hopf algebra over $k$  iff $H^F$ is a
semisimple Hopf algebra over $F$;

(2)  Furthermore, if $H$ is a finite-dimensional Hopf algebra,
then
  $H$ is a cosemisimple  Hopf algebra over $k$  iff $H^F$ is a cosemisimple
Hopf algebra over $F$.

\end {Lemma}
{\bf Proof.} (1) It is clear that $\int _H^l \otimes F = \int
_{H^F}^l.$   Thus
 $H$ is a semisimple  Hopf algebra over $k$  iff $H^F$ is a semisimple
Hopf algebra over $F$.

(2)  $(H \otimes F)^* = H^* \otimes F $ since $H^* \otimes F
\subseteq  (H \otimes F)^*$  and $dim_F (H \otimes F ) = dim_F
(H^* \otimes F)= dim _k H$. Thus we can obtain part (2) by  Part
(1). $\Box$

By the way, if $H$ is a semisimple Hopf algebra, then $H$ is a
separable algebra by Lemma \ref {10.4.11} (see \cite [P284]
{Pa77}).

\begin {Proposition}\label {10.4.12}
Let $F$ be an algebraic closure of $k$, $R$ an algebra over $k$
and
$$r(R \otimes _k F) = r(R) \otimes _k F.$$
 If $H$ is a finite-dimensional  Hopf algebra  with cocommutative
 coradical over $k$ ,
then
     $$r(R) ^{dim H} \subseteq r_H(R).$$

\end {Proposition}
{\bf Proof.} It is clear that $H^F$  is   a finite-dimensional
Hopf algebra
 over $F$ and $dim H = dim  H^F =n$. Let $H^F_0$  be the coradical
 of  $H^F$,  $H^F_1 = H^F_0 \wedge H^F_0,
H^F_{i+1} = H^F_0 \wedge H^F_i$ for $i = 1, \cdots, n-1$. Notice
$H^F_0 \subseteq H_0 \otimes F $. Thus $H^F _0 $  is
cocommutative. It is clear that $H^F_0=kG$  by \cite [Lemma 8.0.1
(c)]{Sw69a} and $H^F= \cup H^F_i$. It is easy to show that if $k >
i $, then  $$H^F_i \cdot (r(R ^F))^k\subseteq r(R ^F)$$ by
induction for $i$. Thus
$$H^F \cdot (r(R ^F))^{ dim H} \subseteq r(R^F),$$
which implies that $(r(R^F))^{dim H }\subseteq r(R^F)_{H^F}$. By
assumption, we have that $(r(R) \otimes F)^{dim H }\subseteq (r(R)
\otimes F)_{H^F}$. It is clear that  $(I \otimes F)_{H^F} = I_H
\otimes F$  for any ideal $I$ of $R$. Consequently, $(r(R))^{dim H
}\subseteq r(R)_H$.   $\Box$

\begin {Theorem}\label {10.4.13}
Let $H$ be a  semisimple, cosemisimple and either commutative or
cocommutative Hopf algebra
 over $k$.
If there exists an algebraic closure $F$ of $k$ such that
$$r_j(R \otimes  F) = r_j(R) \otimes F \hbox { \ \ and \ \ }
r_j((R \# _\sigma H) \otimes  F) = r_j(R \# _\sigma H) \otimes
F,$$ then

(1) $r_j (R) = r_{Hj}(R)= r_{jH}(R); $

 (2)
  $r_{j}(R\# _\sigma  H)= r_{Hj}(R)\# _\sigma H.$
\end {Theorem}
{\bf Proof.}
 (1). By Lemma \ref {10.4.11},
$H^F$ is semisimple and cosemisimple. Considering Corollary \ref
{10.4.9}, we have that
   $r_j (R^F) = r_{H^Fj}(R^F)= r_{jH^F}(R^F).$
On the one hand, by assumption, $r_j(R^F) = r_j(R) \otimes F$. On
the other hand, $r_{jH^F} (R^F) = (r_j (R ) \otimes F)_{H^F} =
r_{jH}(R) \otimes F$. Thus $r_j(R) = r_{jH}(R)$.

(2).  It immediately  follows  from part (1). $\Box$

Considering Theorem \ref {10.4.13} and \cite [Theorem 7.2.13]
{Pa77}, we have

\begin {Corollary}\label {10.4.14}
Let $H$ be a semisimple,  cosemisimple and either commutative or
cocommutative Hopf
 algebra over $k$. If there exists an algebraic closure $F$ of $k$ such that
$F/k$ is separable and algebraic, then

  (1) $r_j (R) = r_{Hj}(R)= r_{jH}(R)$;

  (2)   $r_{j}(R\# _\sigma  H)= r_{Hj}(R)\# _\sigma H.$
                  \end {Corollary}

\begin {Lemma}\label {10.4.15}
(Szasz \cite {Sz82})   $$r_j (R) = r_k(R)$$ holds in the following
three cases:

(1) Every element in $R$ is  algebraic over $k$ (\cite
[Proposition 31.2] {Sz82});

(2) The cardinality of $k$  is strictly greater  than the
dimension of $R$ and $k$ is infinite (\cite [Theorem 31.4]
{Sz82});

(3) $k$ is uncountable and $R$ is  finitely generated (\cite
[Proposition 31.5] {Sz82}).
\end {Lemma}

\begin {Proposition}\label {10.4.16}
Let $F$ be an extension of  $k$. Then $r(R) \otimes  F \subseteq
r(R \otimes  F),$ where $r$ denotes $r_b, r_k, r_l, r_n  .$

\end {Proposition}
{\bf Proof.} When $r =r_n$, for any $x \otimes a  \in r_n (R)
\otimes F$ with $a \not=0$, there exists $y \in R$ such that $x
=xyx$. Thus $x \otimes a = (x \otimes a)(y \otimes a^{-1}) (x
\otimes a)$, which implies $r_n(R) \otimes F \subseteq r_n(R
\otimes F)$.

Similarly, we can obtain the others. $\Box$

\begin {Corollary}\label {10.4.17}
Let $H$ be a semisimple, cosemisimple and commutative or
cocommutative Hopf algebra. If there exists an algebraic closure
$F$ of $k$ such that $F/k$ is a pure transcendental extension and
one of the following three conditions holds:

(i) every element in $R\#_\sigma H$ is  algebraic over $k$;

(ii) the cardinality of $k$  is strictly greater  than the
dimension of $R$ and $k$ is infinite;

(iii) $k$ is uncountable and $R$ is  finitely generated;

then

(1) $r_j (R) = r_{Hj}(R)= r_{jH}(R)$;

(2)    $r_{j}(R\# _\sigma  H)= r_{Hj}(R)\# _\sigma H$;

(3) $r_j(R)= r_k(R) $  and $r_j(R \#_\sigma H) = r_k(R\#_\sigma
H).$
                  \end {Corollary}
{\bf Proof.} First, we have that part (3) holds by Lemma \ref
{10.4.15}. We next see that
\begin {eqnarray*}
r_j(R \otimes F) &\subseteq & r_j(R) \otimes F  \hbox { \ \  \cite
[Theorem
7.3.4] {Pa77}} \\
&=&  r_k(R) \otimes F  \hbox { \ \  part (3)} \\
&\subseteq & r_k(R \otimes F ) \hbox { \ \ proposition \ref {10.4.16} } \\
 &\subseteq & r_j(R \otimes F).
\end {eqnarray*}
Thus $r_j(R \otimes F) = r_j (R) \otimes F.$ Similarly, we can
show that
 $r_j((R \#_\sigma H) \otimes F) = r_j (R\#_\sigma H) \otimes F.$

Finally, using Theorem \ref {10.4.13}, we complete the proof.
$\Box$

\section {The $H$-Von Neumann regular radical}

In this section, we construct the $H$-von Neumann regular radical
for $H$-module algebras and show that it is an $H$-radical
property.

\begin {Definition} \label {10.5.1}
Let $a \in R$. If $a\in (H \cdot a) R (H\cdot a)$, then $a$ is
called an $H$-von Neumann regular element, or an $H$-regular
element in short. If every element of $R$ is an $H$-regular, then
$R$ is called an $H$-regular module algebra, written as
$r_{Hn}$-$H$-module algebra.  $I$ is an $H$-ideal
 of $R$ and every element in $I$ is $H$-regular, then $I$ is called an $H$-
 regular ideal.
\end {Definition}

\begin {Lemma}\label {10.5.2} If $I$ is an $H$-ideal of $R$ and $a \in I$,
then $a$ is $H$-regular in $I$ iff $a$ is $H$-regular in $H$.
\end {Lemma}
{\bf Proof.} The necessity is clear.

Sufficiency: If $a \in (H \cdot a)R (H\cdot a),$ then there exist
$h_i, h_i' \in H, b_i \in R,$  such that
$$a = \sum (h_i \cdot a) b_i (h_i' \cdot a).$$
We see that
\begin {eqnarray*}
a &=& \sum _{i, j} [ h_i \cdot (( h_j \cdot a) b_j (h_j' \cdot
a))]b_i
(h_i' \cdot a) \\
&=&      \sum _{i, j} [ ((h_i)_1 \cdot ( h_j \cdot a)) ((h_i)_2
\cdot b_j)
      ((h_i)_3 \cdot (h_j' \cdot a))]b_i
(h_i' \cdot a)    \\
&\in & (H\cdot a) I (H\cdot a).
\end {eqnarray*}
Thus $a$ is an $H$-regular in $I$. $\Box$

\begin {Lemma}\label {10.5.3}
If  $x - \sum _i (h_i \cdot x) b_i (h_i' \cdot x)$ is $H$-regular,
then $x$ is $H$-regular, where $x, b_i \in R, h_i , h_i' \in H.$
\end {Lemma}
{\bf Proof.} Since   $x - \sum _i (h_i \cdot x) b_i (h_i' \cdot
x)$    is $H$-regular, there exist $g_i, g_i' \in H, c_i \in R $
such that
$$x - \sum _i (h_i \cdot x) b_i (h_i' \cdot x)  =
\sum _j ( g_j \cdot (x - \sum _i (h_i \cdot x) b_i (h_i' \cdot
x))) c_j ( g_j' \cdot ( x - \sum _i (h_i \cdot x) b_i (h_i' \cdot
x))). $$ Consequently, $x \in (H\cdot x)R (H \cdot x).$ $\Box$

\begin {Definition}\label {10.5.4}
$$r_{Hn}(R):=  \{ a \in R \mid \hbox  { \ the \ }   H\hbox {-ideal } (a)
\hbox { generated by } a \hbox { is } H \hbox {-regular } \}.$$
\end {Definition}

\begin {Theorem}\label {10.5.5}
$r_{Hn}(R)$ is an $H$-ideal of $R$.
\end {Theorem}
{\bf Proof.} We first show that $R r_{Hn}(R) \subseteq r_{Hn}(R).$
 For any $a \in r_{Hn}(R), x \in R,$ we have that
 $(xa) $  is $H$-regular since $(xa) \subseteq (a).$
 We next show that $a-b \in r_{Hn}(R)$  for any $a , b \in r_{Hn}(R).$
 For any $x \in (a-b),$ since $(a-b) \subseteq (a) + (b)$, we have that
 $x = u -v$ and $u \in (a), v\in (b).$
 Say $u = \sum _i (h_i \cdot u) c_i (h_i' \cdot u)$  and $h_i, h_i' \in H,
  c_i \in R.$
We see that
\begin {eqnarray*}
&x& - \sum _i (h_i \cdot x)c_i (h_i' \cdot x)\\
 &=& (u-v) - \sum _i  (h_i \cdot (u-v))c_i (h_i' \cdot (u-v))   \\
&=& -v - \sum _i [ -  (h_i \cdot u)c_i (h_i' \cdot v) -   (h_i
\cdot v)c_i (h_i' \cdot u)
+  (h_i \cdot v)c_i (h_i' \cdot v)] \\
& \in & (v) .
\end {eqnarray*}
Thus  $x - \sum _i (h_i \cdot x)c_i (h_i' \cdot x) $ is
$H$-regular and $x$ is $H$-regular by Lemma \ref {10.5.3}.
Therefore $a-b \in r_{Hn}(R).$ Obviously, $r_{Hn}(R)$ is
$H$-stable. Consequently, $r_{Hn}(R)$ is an $H$-ideal of $R$.
$\Box$

\begin {Theorem}\label {10.5.5}
$r_{Hn}(R/r_{Hn}(R)) = 0.$

\end {Theorem}
{\bf Proof.}   Let $\bar R = r/r_{Hn}(R)$ and  $\bar b = b +
r_{Hn}(R) \in r_{Hn}(R/r_{Hn}(R)).$ It is sufficient to show that
$b \in r_{Hn}(R).$  For any $a \in (b)$, it is clear that  $\bar a
\in (\bar b)$.  Thus there exist $h_i, h_i' \in H, \bar c_i \in
\bar R$ such that
$$\bar a = \sum _i (h_i \cdot \bar a) \bar c_i (h_i' \cdot \bar a)
= \sum _i \overline { (h_i \cdot  a)  c_i (h_i' \cdot  a)}.$$ Thus
$ a-  \sum _i  (h_i \cdot  a)  c_i (h_i' \cdot  a) \in r_{Hn}(R),$
which implies that $ a $  is $H$-regular. Consequently,   $b \in
r_{Hn}(R).$ Namely, $\bar b =0$ and $r_{Hn}(R)=0.$  $\Box$

\begin {Corollary}\label {10.5.7.1}
$r_{Hn}$ is an $H$-radical property for $H$-module algebras and
$r_{nH} \le r_{Hn}$.
\end {Corollary}
{\bf Proof.} (R1) If $R \stackrel {f} {\sim }  R'$  and $R$  is an
$r_{Hn}$ -$H$-module algebra, then for any $f(a) \in R'$, $f (a)
\in (H \cdot f(a)) R' (H \cdot f(a)).$  Thus $R'$ is also an
$r_{Hn}$- $H$-module algebra.

(R2) If $I$ is an $r_{Hn}$-$H$-ideal of $R$  and $ r_{Hn}(R)
\subseteq I$ then, for any $a\in I,$  $(a)$  is $H$-regular since
$(a) \subseteq I.$  Thus $I \subseteq r_{Hn}(R).$

  (R3) It follows from Theorem \ref {10.5.6}.

  Consequently $r_{Hn}$ is an $H$-radical property for $H$-module algebras.
  It is straightforward  to check  $r_{nH} \le r_{Hn}.$
 $\Box$

  $r_{Hn}$ is called the  $H$-von Neumann regular radical.

\begin {Theorem}\label {10.5.6}
If $I$ is an $H$-ideal of $R$, then $r_{Hn}(I) = r_{Hn}(R) \cap I
$. Namely, $r_{Hn}$  is a strongly hereditary $H$-radical
property.

\end {Theorem}
{\bf Proof.} By Lemma \ref {10.5.2}, $r_{Hn} (R) \cap I \subseteq
r_{Hn}(I).$ Now, it is sufficient to show that $(x)_I = (x)_R$ for
any $x \in r_{Hn}(I)$, where $(x)_I$ and $(x)_R$ denote the
$H$-ideals generated by $x$ in $I$ and $R$ respectively. Let $x =
\sum (h_i \cdot x) b_i (h_i ' \cdot x)$ , where $h_i, h_i' \in H,
b_i \in I.$  We see that
\begin {eqnarray*}
R(H\cdot x) &=& R( H \cdot ( \sum (h_i \cdot x) b_i (h_i ' \cdot x)) \\
 &\subseteq &R (H\cdot x) I (H \cdot x)  \\
 &\subseteq & I (H \cdot x).
\end {eqnarray*}
Similarly, $$(H\cdot x)R \subseteq (H\cdot x)I.$$ Thus $(x)_I=
(x)_R.$  $\Box$

A graded algebra $R$ of type $G$ is said to be Gr-regular if for
every homogeneous $a \in R_g$ there exists $b\in R$ such that $a=
aba$
 \ \ \ ( see \cite {NO82} P258 ). Now, we give the relations between
 Gr-regularity and  $H$-regularity.

\begin {Theorem}\label {10.5.7}
              If $G$ is a finite group, $R$ is a graded algebra of type $G$,
               and $H=(kG)^*$, then
$R$ is  Gr-regular iff  $R$ is  $H$-regular.

\end {Theorem}
{\bf Proof.} Let $\{ p_g \mid g\in G \}$ be the dual base of base
$\{ g \mid g\in G \}$. If $R$ is Gr-regular for any $a\in R$, then
$a = \sum _{g \in G}  a_g$ with $a_g \in R_g$. Since $R$ is
Gr-regular, there exist $b_{g^{-1}} \in R_{g^{-1}}$ such that $a_g
= a_g b_{g^{-1}} a_g $  and
$$ a = \sum _{g \in G} a_g = \sum _{g \in G } a_g b_{g^{-1}} a_g
= \sum _{g \in G}( p_g \cdot a ) b_{g^{-1}} (p_g \cdot a).$$
Consequently, $R$ is $H$-regular.

Conversely, if $R$ is $H$-regular, then for any $a \in R_g,$ there
exists $b_{x, y}\in R $  such that
$$a = \sum _{x,y \in G} (p_x \cdot a) c_{x,y} (p_y \cdot a).$$
Considering $a \in R_g$, we have that $a = a b_{g,g} a.$ Thus $R$
is Gr-regular .$\Box$

\section { About J.R. Fisher's question }
 In this section,
 we answer the question  J.R. Fisher asked in \cite {Fi75}. Namely,
 we give a necessary and sufficient condition for validity of relation (2) .

Throughout this section, let $k$ be a commutative ring with unit,
$R$ an $H$- module algebra and $H$ a Hopf algebra over $k.$

\begin {Theorem}\label {10.6.1}
Let  ${\cal K}$        be an ordinary special class of rings and
closed with respect to isomorphism. Set $r= r^{\cal K}$ and $\bar
r_H (R) = r(R \# H) \cap R$ for any $H$-module algebra $R$. Then
$\bar r_H$ is an $H$-radical property of $H$-module algebras.
Furthermore, it is an $H$-special radical.

\end {Theorem}

{\bf Proof.} Let $\bar {\cal M}_R = \{ M \mid M $ is an $R$-prime
module and $R/(0:M)_R \in {\cal K } \}$    for any ring $R$ and
$\bar {\cal M} = \cup \bar {\cal M}_R.$ Set ${\cal M}_R =\{ M \mid
M \in \bar {\cal M }_{R \# H} \}$ for any $H$- module algebra $R$
and ${\cal M}= \cup {\cal M}_R.$  It is straightforward to check
that $\bar {\cal M}$  satisfies the conditions of \cite
[Proposition 4.3] {Zh97b}. Thus ${\cal M}$ is an $H$-special
module by      \cite [Proposition 4.3] {Zh97b}.    It is clear
that ${\cal M}(R) = \bar { \cal M} (R \#H) \cap R = r(R \#H) \cap
R$  for any $H$-module algebra $R$. Thus $\bar r_H$ is an
$H$-special radical by \cite [Theorem 3.1] {Zh97b}. $\Box$

Using the Theorem \ref {10.6.1}, we have that
 $\bar r_{bH}, \bar r_{lH}, \bar r_{kH}, \bar r_{jH}, \bar r_{bmH}$
are all $H$-special radicals.

\begin {Proposition}\label {10.6.2}
Let ${\cal K}$ be a  special class of rings and closed with
respect to isomorphism.  Set $r= r^{\cal K}.$ Then

(1) $\bar r_H (R) \# H \subseteq r(R \#H);$

(2) $\bar r_H (R) \# H = r(R \#H)$  iff there exists an $H$-ideal
$I$  of $R$ such that $r(R \#H) = I \#H$;

(3)  $R$ is an $ \bar r_H$-$H$-module algebra iff $r(R \#H) =R\#
H$;

(4)  $I$ is an $ \bar r_H$-$H$-ideal of $R$ iff $r(I \#H)= I \#H$;

(5) $r(\bar r _H (R) \# H) = \bar r _H(R) \# H$ ;

(6) $r(R \# H) = \bar r _H(R) \# H$  iff
 $r(\bar r _H (R) \# H) =  r (R \# H).$

\end {Proposition}

{\bf Proof.} (1)  It is similar to the proof of Proposition \ref
{10.4.3}.

(2) It is a straightforward verification.

(3) If $R$ is an $\bar r_H$-module algebra, then $R \# H \subseteq
r(R \# H)$
 by part (1). Thus $R\# H = r(R \# H)$. The sufficiency is obvious.

(4), (5) and (6) immediately follow from part (3) . $\Box$

\begin {Theorem}\label {10.6.3}  If $R$ is an algebra over field $k$ with
unit and $H$ is a Hopf algebra over field $k$, then

(1)  $\bar r_{jH} (R) = r_{Hj}(R)$ and  $r_j(r_{Hj}(R)\# H) =
r_{Hj}(R) \#H;$

(2)   $r_j(R\# H) = r_{Hj}(R) \#H$  iff $r_j(r_{Hj}(R)\# H) =
r_{j}(R \#H)$    iff $r_j(r_j(R \# H)\cap R \# H) = r_j(R\# H);$

(3) Furthermore, if   $H$ is finite-dimensional, then
    $r_j(R\# H) = r_{Hj}(R) \#H$  { \ \ \ \ } iff \\
$r_j(r_{jH}(R)\# H) = r_{j}(R \#H).$

\end {Theorem}

{\bf Proof.}  (1) By \cite [Proposition 3.1 ]  {Zh98a}, we have
$\bar r_{jH} (R) = r_{Hj}(R)$. Consequently, $r_j(r_{Hj}(R)  \# H)
=r_{Hj}(R) \#H$ by Proposition \ref {10.6.2} (5).

(2) It immediately  follows from part (1) and Proposition \ref
{10.6.2} (6).

(3) It can easily be proved by part (2) and \cite [Proposition 3.3
(2)] {Zh98a}. $\Box$

The theorem answers the question J.R. Fisher asked in \cite {Fi75}
: When is $r_j(R \#H) = r_{Hj}(R)\#H$ ?

\begin {Proposition}\label {10.6.4}  If $R$ is an algebra over field $k$ with
unit and $H$ is a finite-dimensional Hopf algebra over field $k$,
then

(1)  $\bar r_{bH} (R) = r_{Hb}(R)=r_{bH}(R)$ and  $r_b(r_{Hb}(R)\#
H) = r_{Hb}(R) \#H;$

(2)   $r_b(R\# H) = r_{Hb}(R) \#H$  iff $r_b(r_{Hb}(R)\# H) =
r_{b}(R \#H)$  iff $r_b(r_{bH}(R)\# H) = r_{b}(R \#H)$  iff
$r_b(r_{b}(R \#H) \cap R \# H) = r_{b}(R \#H).$

\end {Proposition}

{\bf Proof.}  (1) By \cite [Proposition 2.4 ]  {Zh98a}, we have
$\bar r_{bH} (R) = r_{Hb}(R)$.  Thus $r_b(r_{Hb}(R)\# H) =
r_{Hb}(R )\#H$  by Proposition \ref {10.6.2} (5).

(2) It follows from part (1) and Proposition \ref {10.6.2} (6) .
$\Box$

In fact,  if $H$ is commutative or cocommutative, then $S^2 =
id_H$
 by \cite [Proposition 4.0.1] {Sw69a},
and  $H$  is semisimple and cosemisimple iff the character $char k
$ of $k$ does  not divides $dim H$   \ \ ( see \cite [Proposition
2 (c)] {Ra94} ). It is clear that if $H$ is a finite-dimensional
commutative or cocommutative Hopf algebra and
 the character $char k $ of $k$ does  not
divides $dim H$, then  $H$ is a finite-dimensional semisimple,
cosemisimple, commutative or cocommutative Hopf algebra.
Consequently, the conditions in Corollary \ref {10.4.9}, Theorem
\ref {10.4.13}, Corollary  \ref{10.4.14} and \ref {10.4.17}
  can be simplified

\end {document}